\documentclass[12pt]{article}

\usepackage{amsmath}
\numberwithin{equation}{section}
\allowdisplaybreaks
\usepackage{amsfonts}
\usepackage{amssymb}
\usepackage{enumitem}
\usepackage[colorlinks=true,linkcolor=black,citecolor=blue,urlcolor=blue,]{hyperref}
\newtheorem{theorem}{Theorem}[section]
\newtheorem{lemma}{Lemma}[section]
\newtheorem{corollary}{Corollary}[section]

\newenvironment{proof}{\noindent{\bf Proof:}}{\hfill\fbox{}\vspace*{1mm}}

\begin{document}
	
	\title{Difference of Composition Operators on Korenblum Spaces over Tubular Domains}

	\author{ Yuheng Liang, Lvchang Li, Haichou Li*}
	\author{Yuheng Liang \thanks
		{College of Mathematics and Informatics,
			South China Agricultural University,
			Guangzhou,
			510640,
			China
			Email: 20232115004@stu.scau.edu.cn.},\
		Lvchang Li\thanks{College of Mathematics and Informatics,
			South China Agricultural University,
			Guangzhou,
			510640,
			China
			Email: 20222115006@stu.scau.edu.cn.},\
		Haichou Li* \thanks{ Corresponding author,
			College of Mathematics and Informatics,
			South China Agricultural University,
			Guangzhou,
			510640,
			China
			Email: hcl2016@scau.edu.cn. }
	}
	
	\date{}
	\maketitle
	\begin{center}
		\begin{minipage}{120mm}
			\begin{center}{\bf Abstract}\end{center}
			{This paper investigates the properties of the difference of composition operators on the Korenblum space over the product of upper half planes, characterizing their boundedness and compactness. Using the result on boundedness, we show that all bounded differences of composition operators are absolutely summing operators.}		
			
		\end{minipage}
	\end{center}
	
	\maketitle

\section{Introduction}
\ \ \ \
The study of composition operators in analytic function space, particularly on spaces over the unit disc $\mathbb{D}$, has a long-standing history and plays a crucial role in addressing various problems in fields such as functional analysis, complex dynamical systems, representation theory, number theory, and mathematical physics. Notably, J. H. Shapiro and other mathematicians have made significant contributions to problems related to composition operators, yielding many valuable results (see for example \cite{CCOBG}, \cite{BO HO}, \cite{GB upper}, \cite{CO Wa Ds} etc. for some interesting results on composition operators). Over the past forty years, research has increasingly focused on the topological structures of the sets of composition operators. To investigate these topological structures, understanding the properties of the difference of composition operators has become essential (see such as \cite{topo DCO}). The boundedness, compactness, and many other properties of differences of composition operators are now central topics of interest among operator theorists, see for example \cite{DCO disc}.

Apart from the unit disc, current studies in operator theory also explores function spaces over more complex domains. One generalization is from bounded domains to unbounded domains, such as from the unit disc $\mathbb{D}$ to the upper half-plane $\mathcal{H}$. In \cite{BHW}, B. R. Choe et al. consider the difference of composition operators on weight Bergman spaces over upper half-plane, have obtained many valuable results. C. B. Pang and M. F. Wang characterized in \cite{DCO Half} the boundedness and compactness of the difference of composition operators from weighted Bergman space to Lebesgue space over upper half-plane. Several years later, Wang and X. Guo discussed in \cite{main} the related problems on Korenblum spaces. Besides upper half-plane, the right half-plane is often studied too (see \cite{CSDWCO}, \cite{CO Wa Ds}, etc.) and it can be associated with many problems in number theory. Another generalization is from $1$-dimensional spaces to multidimensional spaces. In bounded cases, one may consider the ball $\mathbb{B}^n$ or polydisc $\mathbb{D}^n$ (see such as \cite{DWCO polydisk}), while in unbounded cases, tube domains are studied. The domain $\mathcal{H}^n$ considered in this paper is a special type of tube domain in $\mathbb{C}^n$. 

The Korenblum space, often referred to as a growth space, is a special type of analytic function space. It is well-known that Korenblum space is a Banach space, and it plays a significant role in the study of the interpolation problem in Bergman spaces.

Another class of operators that has also been widely studid is the absolutely summing operators. Since Grothendieck, research on absolutely summing operators has been ongoing. Later, the concept of absolute summing operators was extended to the concept of $p$-summing operators by Pietsch, and many intriguing properties of these classes of operators had been revealed. A standard reference in this area is \cite{ASO}. The $p$-summing operators have significant applications in the study of operator theory and the geometry of Banach spaces, see \cite{ASO}, \cite{ASO Lp} or \cite{BSA}. T. Domenig studied in \cite{Domenig} the $p$-summing composition operators on Bergman spaces. Recently, B. He et al. gave from another point of view a new approach in \cite{ASCE} to the characterization of $p$-summing weighted composition operators on Bergman spaces.

Motivated by the aforementioned studies, especially by \cite{main}, we attempt to generalize the relevant results to higher-dimensional spaces. This article primarily explores the boundedness, compactness, and $p$-summing property.

This paper is organized as follows. In section 2 we introduce some fundamental concepts and useful tools with their basic properties. Some lemmas as preparations for our proof of the main results are provided in section 3. In section 4, we discuss the necessary and sufficient conditions for when the difference of composition operators is bounded from $\mathcal{A}^{-\gamma}(\mathcal{H}^n)$ to $\mathcal{A}^{-\gamma}(\mathcal{H}^n)$. As an application, we prove in the same section that every bounded difference of composition operators from $\mathcal{A}^{-\gamma}(\mathcal{H}^n)$ into $\mathcal{A}^{-\gamma}(\mathcal{H}^n)$ is an absolutely summing operator. In the last section, we characterize the compactness of these operators.

\section{Preliminaries}
\ \ \ \
Let $\mathcal{H}$ be the upper half-plane in $\mathbb{C}$ and by
$$\mathcal{H}^n:=\{z=(z_1,\cdots,z_n)\in\mathbb{C}^n:\mathrm{Im} z_k>0,k=1,\cdots,n\}$$
we denote the Cartesian product of $n$ upper half planes in $\mathbb{C}^n$, where $\mathbb{C}^n$ is the $n$-dimensional Euclidean complex space. This is one of the natural generalizations of the upper half-plane to higher-dimensional spaces, and it is a tube domain in $\mathbb{C}^n$. Further results on $\mathcal{H}^n$ can be found in such as \cite{multihalfplane}.

The set $H(\mathcal{H}^n)$ is the collection of all holomorphic functions on $\mathcal{H}^n$ and $\mathcal{S}(\mathcal{H}^n)$ the collection of all holomorphic self-maps of $\mathcal{H}^n$. For any $\varphi\in\mathcal{S}(\mathcal{H}^n)$, we define the composition operator $C_{\varphi}$ by
$$(C_{\varphi} f)(z)=f(\varphi(z))$$
for $f\in H(\mathcal{H}^n)$. Suppose that $\varphi=(\varphi_1,\cdots,\varphi_n), \psi=(\psi_1,\cdots,\psi_n)\in\mathcal{S}(\mathcal{H}^n)$, we consider the difference of $C_{\varphi}$ and $C_{\psi}$ in this paper .

If $\gamma=(\gamma_1,\cdots,\gamma_n)$ is a multi-index, then we denote
\[|\gamma|:=\gamma_1+\cdots+\gamma_n.\]

For any $\gamma\in\mathbb{R}_+^n=\{\gamma=(\gamma_1,\cdots,\gamma_n)\in\mathbb{R}^n:\gamma_k>0,k=1,\cdots,n\}$, we define the Korenblum space over $\mathcal{H}^n$ to be
$$\mathcal{A}^{-\gamma}(\mathcal{H}^n):=\{f\in H(\mathcal{H}^n):\lVert f \rVert_{\mathcal{A}^{-\gamma}(\mathcal{H}^n)}=\sup_{z\in \mathcal{H}^n}\left|f(z)\right|\prod_{k=1}^n (\mathrm{Im} z_k)^{\gamma_k}<\infty\}.$$
Here the vector-type parameter $\gamma$ is a multi-index.

Recall that, for $z,w\in\mathcal{H}$, we frequently consider the pseudo-hyperbolic distance $d(z,w)$ between them, which is given by
$$d(z,w):=\left|\frac{z-w}{z-\overline{w}}\right|,$$
and a straightforward calculation shows that
\begin{equation}\label{rho ineq1}
	\frac{1-\delta}{1+\delta}<\frac{\mathrm{Im} z}{\mathrm{Im} w}<\frac{1+\delta}{1-\delta}
\end{equation}
and
\begin{equation}\label{rho ineq2}
	\frac{1-\delta}{1+\delta}<\frac{\left|z-\overline{w}\right|}{2\mathrm{Im} z}<\frac{1+\delta}{1-\delta}
\end{equation}
whenever $d(z,w)<\delta$. Moreover, we have
\begin{equation}\label{rho ineq3}
	\frac{1-\delta}{1+\delta}<\left|\frac{z-\overline{a}}{w-\overline{a}}\right|<\frac{1+\delta}{1-\delta}
\end{equation}
whenever $d(z,w)<\delta$ and $a\in\mathcal{H}$. See \cite{pho ineq} for more details.

For $z=(z_1,\cdots,z_n)$, $w=(w_1,\cdots,w_n)\in\mathcal{H}^n$, we write $\rho_k(z,w):=d(z_k,w_k)$. The pseudo-hyperbolic distance in $\mathcal{H}^n$ is defined by
$$\rho(z,w):=\max_{1\leq k\leq n}\rho_k(z,w).$$
For convenience, we denote $\rho(z):=\rho(\varphi(z),\psi(z))$ and $\rho_k(z):=\rho_k(\varphi(z),\psi(z))$.

For $\delta=(\delta_1,\cdots,\delta_n)$ where $0<\delta_k<1,k=1,\cdots,n$ and $z\in\mathcal{H}^n$, the pseudohyperbolic polydisc centered at $z$ with multi-radius $\delta$ is given by
$$E_{\delta}(z):=\{w\in\mathcal{H}^n:\rho_k(z,w)<\delta_k,k=1,\cdots,n\}.$$
In fact, $E_{\delta}(z)$ is the Euclidean polydisc 
$$\left\{w\in\mathbb{C}^n:\left|w_k-(\mathrm{Re} z_k+i\frac{1+\delta_k^2}{1-\delta_k^2}\mathrm{Im} z_k)\right|<\frac{2\delta_k}{1-\delta_k^2}\mathrm{Im} z_k,k=1,\cdots,n\right\}$$
in $\mathbb{C}^n$, and its boundary
$$\partial E_{\delta}(z):=\{w\in\mathcal{H}^n:\rho_k(z,w)=\delta_k,k=1,\cdots,n\}$$
is an $n$-dimensional Euclidean torus.

We denote  by $z\to\partial\widehat{\mathcal{H}^n}$ the statement that $z_k\to0$ or $z_k\to\infty$ for some $k=1,2,\cdots,n$. When we say $\lim_{z\to\partial\widehat{\mathcal{H}^n}}f(z)=0$, it equivalently means that for any $\varepsilon>0$, there exists a compact subset $\Omega\subset\mathcal{H}^n$ such that $\sup_{z\in\mathcal{H}^n\setminus\Omega}|f(z)|<\varepsilon$.

At the end of this section, we briefly introduce the concept of absolutely summable operators. Suppose that $1\leq p<\infty$ and let $X$, $Y$ be Banach spaces, we say a linear operator $T:X\to Y$ is $p$-summing if there is a constant $C_p>0$ such that for all $m\in\mathbb{N}$ and $x^1,\cdots,x^m\in X$ we have
$$\left(\sum_{j=1}^{m}\|Tx^j\|^p\right)^{\frac{1}{p}}\leq C_p \sup_{u\in B_{X^*}}\left(\sum_{j=1}^{m}|u(x^j)|^p\right)^{\frac{1}{p}},$$
where $B_{X^*}$ is the closed unit ball in the dual space $X^*$ of $X$. The least choice of $C_p$ is always denoted by $\pi_p(T)$. The case of $1$-summing is also called absolutely summing and we sometimes write $\pi(T):=\pi_p(T)$. In other words, an operator from $X$ to $Y$ is absolutely summing if it takes unconditionally summable sequences $\{x_j\}$ in $X$ to absolutely summable sequences $\{Tx_j\}$ in $Y$.

Here are some basic facts of $p$-summing\textsuperscript{\cite{ASO}}.

\begin{enumerate}
	\item[(a)] The set $\Pi_p(X,Y)$ of all $p$-summing operators from $X$ to $Y$ is a Banach subspace of the space $B(X,Y)$ of all continuous linear operators. Moreover, $\pi_p(\cdot)$ defines a norm on $\Pi_p(X,Y)$ with 
	$$\|T\|_{B(X,Y)}\leqslant\pi_p(T),\quad T\in\Pi_p(X,Y).$$
	\item[(b)] If $1\leqslant p\leqslant q<\infty$, then $\Pi_p(X,Y)\subset\Pi_q(X,Y)$ and for any $T\in\Pi_q(X,Y)$, we have
	$$\pi_q(T)\leqslant\pi_p(T).$$
\end{enumerate}

From this we can see that if $T$ is $p$-summing for some $p\geqslant 1$ then $T$ must be a bounded operator, and that if $T$ is absolutely summing then $T$ is $p$-summing for any $1<p<\infty$.

Throughout the paper, the notation $A\lesssim B$ means that there is a positive constant $C$ such that $A\le CB$, and the notation $A\simeq B$ means that $A\lesssim B$ and $B\lesssim A$.

\section{Basic lemmas}
\ \ \ \ 
In this section, we will introduce several key lemmas, which will play an important role in the proof of the main theorems.

The first lemma generalize the lemma 7.2 in \cite{BHW}.

\begin{lemma}\label{mi yu 1 cha}
	Given $s,\delta\in\mathbb{R}_+^n$ where $0<\delta_k<1$ for any $k=1,\cdots,n$, there is a constant $C_{s,\delta}>0$ such that
	$$\left|\prod_{k=1}^n \left(\frac{z_k-\overline{a_k}}{w_k-\overline{a_k}}\right)^{s_k}-1\right|\leq C_{s,\delta}\rho(z,w)$$
	for all $a,z,w\in\mathcal{H}^n$ with $\rho_k(z,w)<\delta_k$.
\end{lemma}

\begin{proof}
	By (\ref{rho ineq3}) we have
	$$\frac{1-\delta_k}{1+\delta_k}<\left|\frac{z_k-\overline{a_k}}{w_k-\overline{a_k}}\right|<\frac{1+\delta_k}{1-\delta_k}.$$
	
	Consider the function on $\mathcal{H}^n$ defined by
	\[h(z)=
	\begin{cases}
		\displaystyle \frac{\prod_{k=1}^{n} z_k^{s_k}-1}{\max_{1\leq k\leq n}|z_k-1|}\quad&,z\neq (1,\cdots,1) \\
		1&,z=(1,\cdots,1)
	\end{cases}
	\]
	which is bounded on any bounded domain in $\mathbb{C}^n$ hence

	$$M:=\sup\left\{\left|h(\xi)\right|:\frac{1-\delta_k}{1+\delta_k}<\left|\xi_k\right|<\frac{1+\delta_k}{1-\delta_k},\xi_k\notin\mathbb{R}_-,k=1,\cdots,n\right\}<\infty$$
	Then, by (\ref{rho ineq3}), we have
	\begin{align*}
		\left|\prod_{k=1}^{n}\left(\frac{z_k-\overline{a_k}}{w_k-\overline{a_k}}\right)^{s_k}-1\right|=&\left|\frac{\prod_{k=1}^{n}\left(\frac{z_k-\overline{a_k}}{w_k-\overline{a_k}}\right)^{s_k}-1}{\max_{1\leq k\leq n}\left|\frac{z_k-\overline{a_k}}{w_k-\overline{a_k}}-1\right|}\right|\max_{1\leq k\leq n}\left|\frac{z_k-w_k}{w_k-\overline{a_k}}\right|\\	
		\leq&M\max_{1\leq k\leq n}\left|\frac{z_k-w_k}{w_k-\overline{a_k}}\right|\\
		=&M\max_{1\leq k\leq n}\left|\frac{z_k-w_k}{z_k-\overline{w_k}}\cdot\frac{z_k-\overline{w_k}}{w_k-\overline{a_k}}\right|\\
		\leq&M\rho(z,w)\max_{1\leq k\leq n}\left|\frac{z_k-\overline{w_k}}{\mathrm{Im} w_k}\right|\\
		\lesssim &\rho(z,w),
	\end{align*}
where the penultimate step relies on the fact that $\left|w_k-\overline{a}_k\right|>\mathrm{Im} w_k$ and the final step follows from (\ref{rho ineq2}).
\end{proof}	

A practical modification of Lemma \ref{mi yu 1 cha} is the following lemma.

\begin{lemma}\label{xubushang mi 1 cha}
	Given $s,\delta\in\mathbb{R}_+^n$ where $0<\delta_k<1$ for any $k=1,\cdots,n$, there is a constant $C_{s,\delta}>0$ such that
	$$\left|\prod_{k=1}^n \left(\frac{\mathrm{Im} z_k}{\mathrm{Im} w_k}\right)^{s_k}-1\right|\leq C_{s,\delta}\rho(z,w)$$
	for all $z,w\in\mathcal{H}^n$ with $\rho_k(z,w)<\delta_k$.
\end{lemma}

\begin{proof}
	By Lemma \ref{mi yu 1 cha}, we have
	$$
	\begin{aligned}
		&\left|\prod_{k=1}^n \left(\frac{\mathrm{Im} z_k}{\mathrm{Im} w_k}\right)^{s_k}-1\right|\\
		\leq&\left|\prod_{k=1}^n \left(\frac{\mathrm{Im} z_k}{\mathrm{Im} w_k}\right)^{s_k}-\prod_{k=1}^n \left(\frac{z_k-\overline{z_k}}{w_k-\overline{z_k}}\right)^{s_k}\right|+\left|\prod_{k=1}^n \left(\frac{z_k-\overline{z_k}}{w_k-\overline{z_k}}\right)^{s_k}-1\right|\\
		\leq&\prod_{k=1}^n \left(\frac{\mathrm{Im} z_k}{\mathrm{Im} w_k}\right)^{s_k}\left|1- \prod_{k=1}^n\left(\frac{w_k- \overline{w_k}}{w_k-\overline{z_k}}\right)^{s_k}\right|+\left|\prod_{k=1}^n \left(\frac{z_k-\overline{z_k}}{w_k-\overline{z_k}}\right)^{s_k}-1\right|\\
		\lesssim&\prod_{k=1}^n \left(\frac{1+\delta_k}{1-\delta_k}\right)^{s_k}\left|1- \prod_{k=1}^n\left(\frac{w_k- \overline{w_k}}{w_k-\overline{z_k}}\right)^{s_k}\right|+\left|\prod_{k=1}^n \left(\frac{z_k-\overline{z_k}}{w_k-\overline{z_k}}\right)^{s_k}-1\right|\\
		\lesssim& \rho(z,w).
	\end{aligned}
	$$
\end{proof}

From lemma \ref{xubushang mi 1 cha}, it's easy to obtain the next lemma. This inequality is frequently used in this article.

\begin{lemma}\label{micha youxian}
	If $z\in\mathcal{H}^n$, $\varphi, \psi\in\mathcal{S}(\mathcal{H}^n)$, $\gamma\in\mathbb{R}_+^n$, then 
	$$
	\begin{aligned}
	&\left|\prod_{k=1}^n (\frac{\mathrm{Im} z_k}{\mathrm{Im}\varphi_k(z)})^{\gamma_k}-\prod_{k=1}^n (\frac{\mathrm{Im} z_k}{\mathrm{Im}\psi_k(z)})^{\gamma_k}\right|\\
	&\lesssim \rho(z)\left(\prod_{k=1}^n \left(\frac{\mathrm{Im} z_k}{\mathrm{Im} \varphi_k(z)}\right)^{\gamma_k}+\prod_{k=1}^n \left(\frac{\mathrm{Im} z_k}{\mathrm{Im} \psi_k(z)}\right)^{\gamma_k}\right).
	\end{aligned}
	$$
\end{lemma}

Let $\Omega$ be a subset of $\mathcal{H}^n$,
$$S_{\Omega,f}^{\gamma}:=\sup_{z\in\Omega}|f(z)|\prod_{k=1}^n (\mathrm{Im} z_k)^{\gamma_k}.$$
Obviously we have $S_{\Omega_1,f}^{\gamma}\leq S_{\Omega_2,f}^{\gamma}$ for $\Omega_1\subseteq\Omega_2$ and $S_{\mathcal{H}^n,f}^{\gamma}=\|f\|_{\mathcal{A}^{-\gamma}(\mathcal{H}^n)}$. By adapting the proof of Lemma 5.1 in \cite{Bergman zhu}, we can prove the following result.

\begin{lemma}\label{jiaquan hanshu mi cha}
	Let $\gamma\in\mathbb{R}_+^n$, $\Omega$ be an open subset in $\mathcal{H}^n$, then there is a constant $C_{\gamma,\Omega}$ such that 
	$$\left|f(z)\prod_{k=1}^n (\mathrm{Im} z_k)^{\gamma_k}-f(w)\prod_{k=1}^n (\mathrm{Im} w_k)^{\gamma_k}\right|\leq C_{\gamma,\Omega}S_{\Omega,f}^{\gamma}\rho(z,w)$$
	for all $z,w\in\Omega$ and $f\in\mathcal{A}^{-\gamma}(\mathcal{H}^n)$.
\end{lemma}

\begin{proof}

	Let $z\in\Omega$, and $E_{\delta}(z)\subset\Omega$, as given in section 2,  be an open pseudohyperbolic polydisc centered at $z$ with significantly small multi-radius $\delta$ such that $\partial E_{\delta}(z)\subset\Omega$, then $\left|f(\zeta)\right|\leq S_{\Omega,f}^{\gamma}\prod_{j=1}^{n}(\mathrm{Im} \zeta_j)^{-\gamma_j}\sim S_{\Omega,f}^{\gamma}\prod_{j=1}^{n}(\mathrm{Im} z_j)^{-\gamma_j}$ whenever $\zeta\in\overline{E_{\delta}(z)}$, and by Cauchy's formula, we have
	$$\frac{\partial}{\partial z_k}f(z)=\frac{1}{(2\pi i)^n}\int_{\partial E_{\delta}(z)}\frac{f(\zeta)\:d\zeta}{(\zeta_k-z_k)\prod_{j=1}^{n}(\zeta_j-z_j)},\quad z\in\Omega,$$
	
	which yields the estimate
	$$
	\begin{aligned}
		\left|\frac{\partial}{\partial z_k}f(z)\right|\leq&\frac{1}{(2\pi)^n}\int_{\partial E_{\delta}(z)}\frac{\left|f(\zeta)\right|\left|\:d\zeta\right|}{\left|\zeta_k-z_k\right|\prod_{j=1}^{n}\left|\zeta_j-z_j\right|}\\
		\lesssim&\int_{\partial E_{\delta}(z)}\frac{S_{\Omega,f}^{\gamma}\prod_{j=1}^{n}(\mathrm{Im} z_j)^{-\gamma_j}\left|\:d\zeta\right|}{(\mathrm{Im} z_k)\prod_{j=1}^{n}\mathrm{Im} z_j}\\
		=&S_{\Omega,f}^{\gamma}(\mathrm{Im} z_k)^{-1}\prod_{j=1}^{n}(\mathrm{Im} z_j)^{-\gamma_j-1}\int_{\partial E_{\delta}(z)}\left|\:d\zeta\right| \\
		\lesssim& S_{\Omega,f}^{\gamma}(\mathrm{Im} z_k)^{-1}\prod_{j=1}^{n}(\mathrm{Im} z_j)^{-\gamma_j}.
	\end{aligned}
	$$
	If $\rho_k(z,w)$ is sufficiently small, then
	$$|z_k-w_k|\sim\rho_k(z,w)\mathrm{Im} z_k.$$
	
	Now consider $d\left(f(\eta)\prod_{k=1}^n (\mathrm{Im} \eta_k)^{\gamma_k}\right)$ which is the total differential of $f(\eta)\prod_{k=1}^n (\mathrm{Im} \eta_k)^{\gamma_k}$ at $\eta=(\eta_1,\cdots,\eta_n)\in\Omega$, then
	$$
	\begin{aligned}
	&\left|d\left(f(\eta)\prod_{k=1}^n (\mathrm{Im} \eta_k)^{\gamma_k}\right)\right|\\
	&\leqslant\sum_{k=1}^{n}\left|\frac{\partial}{\partial \eta_k}\left(f(\eta)\prod_{k=1}^n (\mathrm{Im} \eta_k)^{\gamma_k}\right)\right|\left|d\eta_k\right| \\
	&\leqslant\sum_{k=1}^{n}\left(\frac{\gamma_k}{2}(\mathrm{Im} \eta_k)^{-1}|f(\eta)|\prod_{j=1}^{n}(\mathrm{Im} \eta_j)^{\gamma_j}+\left|\frac{\partial}{\partial \eta_k}f(\eta)\right|\prod_{j=1}^{n}(\mathrm{Im} \eta_j)^{\gamma_j}\right)|d\eta_k| \\
	&\lesssim S_{\Omega,f}^{\gamma}\sum_{k=1}^{n}(\mathrm{Im} \eta_k)^{-1}|d\eta_k|.
	\end{aligned}
	$$
	It follows that
	$$
	\begin{aligned}
	\left|f(z)\prod_{k=1}^n (\mathrm{Im} z_k)^{\gamma_k}-f(w)\prod_{k=1}^n (\mathrm{Im} w_k)^{\gamma_k}\right|&\lesssim S_{\Omega,f}^{\gamma}\sum_{k=1}^{n}(\mathrm{Im} z_k)^{-1}|z_k-w_k| \\
	&\lesssim S_{\Omega,f}^{\gamma}\sum_{k=1}^{n}\rho_k(z,w) \\
	&\lesssim S_{\Omega,f}^{\gamma}\rho(z,w).
	\end{aligned}
	$$
	This completes the proof.
\end{proof}

When we choose $\Omega=\mathcal{H}^n$, the lemma above is a higher-dimensional analogy of the lemma 2.5 in \cite{main}.

Suppose $a\in\mathcal{H}^n$, $m\in\{1,\cdots,n\}$, two classes of functions below play an important role in this article:
\begin{equation}\label{f_a}
	f_a(z)=\prod_{k=1}^n \frac{((2i)^2\mathrm{Im} a_k)^{\gamma_k}}{(z_k-\overline{a_k})^{2\gamma_k}},z\in\mathcal{H}^n,
\end{equation}
\begin{equation}\label{g_a,m}
	g_{a,m}(z)=4^{|\gamma|}\cdot\frac{z_m-a_m}{z_m-\overline{a_m}}\cdot\prod_{k=1}^n \frac{(\mathrm{Im} a_k)^{\gamma_k}}{(z_k-\overline{a_k})^{2\gamma_k}},z\in\mathcal{H}^n.
\end{equation}

It's easy to show that $f_a,g_{a,m}\in\mathcal{A}^{-\gamma}(\mathcal{H}^n)$ with
\begin{equation}
	\|f_a\|_{\mathcal{A}^{-\gamma}(\mathcal{H}^n)}=|f_a(a)|\prod_{k=1}^{n}(\mathrm{Im} a_k)^{\gamma_k}=1,
\end{equation}
\begin{equation}
	\sup_{a\in \mathcal{H}^n}\|g_{a,m}\|_{\mathcal{A}^{-\gamma}(\mathcal{H}^n)}<\infty
\end{equation}
Besides, $\{f_a\}_a$ and $\{g_{a,m}\}_a$ exhibit nice convergency, namely:
\begin{lemma}
	Let $\Omega$ be any compact subset of $\mathcal{H}^n$, then $f_a$ and $g_{a,m}$ converge uniformly to $0$ on $\Omega$ when $a\to \partial\widehat{\mathcal{H}^n}$.
\end{lemma}

To prove the compactness of $C_{\varphi}-C_{\psi}$ on $\mathcal{A}^{-\gamma}(\mathcal{H}^n)$, the following lemma is useful.

\begin{lemma}\label{limit}
Suppose $z\in\mathcal{H}^n$, $\varphi, \psi\in\mathcal{S}(\mathcal{H}^n)$, $\gamma\in\mathbb{R}_+^n$ and $\Omega_1,\Omega_2$ are open subsets of $\mathcal{H}^n$ such that $\varphi(z)\in\Omega_1$ and $\psi(z)\in\Omega_2$, then
\begin{align*}
	&\left|f(\varphi(z))-f(\psi(z))\right|\prod_{k=1}^n (\mathrm{Im} z_k)^{\gamma_k} \\
	\lesssim& \left(h_1 S_{\Omega_2,f}^{\gamma}+h_2 S_{\Omega_2,f}^{\gamma}\right)\rho(z)\left(\prod_{k=1}^n \left(\frac{\mathrm{Im} z_k}{\mathrm{Im} \varphi_k(z)}\right)^{\gamma_k}+\prod_{k=1}^n \left(\frac{\mathrm{Im} z_k}{\mathrm{Im} \psi_k(z)}\right)^{\gamma_k}\right) \\
	&+S_{\Omega_1\cup\Omega_2,f}^{\gamma}\rho(z)\left(h_2\prod_{k=1}^n \left(\frac{\mathrm{Im} z_k}{\mathrm{Im} \varphi_k(z)}\right)^{\gamma_k}+h_1\prod_{k=1}^n \left(\frac{\mathrm{Im} z_k}{\mathrm{Im} \psi_k(z)}\right)^{\gamma_k}\right)
\end{align*}
for any $f\in\mathcal{A}^{-\gamma}(\mathcal{H}^n)$, $h_1,h_2\geq 0$, $h_1+h_2=1$.
\end{lemma}

\begin{proof}
	A straightforward verification shows that
	\begin{align*}
		&\left|f(\varphi(z))-f(\psi(z))\right|\prod_{k=1}^n (\mathrm{Im} z_k)^{\gamma_k}\\
		=&\left|h_1f(\varphi(z))\prod_{j=1}^n (\mathrm{Im} \varphi_j(z))^{\gamma_j}\left(\prod_{k=1}^n \left(\frac{\mathrm{Im} z_k}{\mathrm{Im} \varphi_k(z)}\right)^{\gamma_k}-\prod_{k=1}^n \left(\frac{\mathrm{Im} z_k}{\mathrm{Im} \psi_k(z)}\right)^{\gamma_k}\right)\right. \\
		&+h_1\prod_{j=1}^n \left(\frac{\mathrm{Im} z_j}{\mathrm{Im} \psi_j(z)}\right)^{\gamma_j}\left(f(\varphi(z))\prod_{k=1}^n (\mathrm{Im} \varphi_k(z))^{\gamma_k}-f(\psi(z))\prod_{k=1}^n (\mathrm{Im} \psi_k(z))^{\gamma_k}\right) \\
		&+h_2\left|f(\psi(z))\right|\prod_{j=1}^n (\mathrm{Im} \psi_j(z))^{\gamma_j}\left(\prod_{k=1}^n \left(\frac{\mathrm{Im} z_k}{\mathrm{Im} \varphi_k(z)}\right)^{\gamma_k}-\prod_{k=1}^n \left(\frac{\mathrm{Im} z_k}{\mathrm{Im} \psi_k(z)}\right)^{\gamma_k}\right) \\
		&+h_2\prod_{j=1}^n \left.\left(\frac{\mathrm{Im} z_j}{\mathrm{Im} \varphi_j(z)}\right)^{\gamma_j}\left(f(\varphi(z))\prod_{k=1}^n (\mathrm{Im} \varphi_k(z))^{\gamma_k}-f(\psi(z))\prod_{k=1}^n (\mathrm{Im} \psi_k(z))^{\gamma_k}\right)\right| \\
		\leqslant& h_1\left|f(\varphi(z))\right|\prod_{j=1}^n (\mathrm{Im} \varphi_j(z))^{\gamma_j}\left|\prod_{k=1}^n \left(\frac{\mathrm{Im} z_k}{\mathrm{Im} \varphi_k(z)}\right)^{\gamma_k}-\prod_{k=1}^n \left(\frac{\mathrm{Im} z_k}{\mathrm{Im} \psi_k(z)}\right)^{\gamma_k}\right| \\
		&+h_1\prod_{j=1}^n \left(\frac{\mathrm{Im} z_j}{\mathrm{Im} \psi_j(z)}\right)^{\gamma_j}\left|f(\varphi(z))\prod_{k=1}^n (\mathrm{Im} \varphi_k(z))^{\gamma_k}-f(\psi(z))\prod_{k=1}^n (\mathrm{Im} \psi_k(z))^{\gamma_k}\right| \\
		&+h_2\left|f(\psi(z))\right|\prod_{j=1}^n (\mathrm{Im} \psi_j(z))^{\gamma_j}\left|\prod_{k=1}^n \left(\frac{\mathrm{Im} z_k}{\mathrm{Im} \varphi_k(z)}\right)^{\gamma_k}-\prod_{k=1}^n \left(\frac{\mathrm{Im} z_k}{\mathrm{Im} \psi_k(z)}\right)^{\gamma_k}\right| \\
		&+h_2\prod_{j=1}^n \left(\frac{\mathrm{Im} z_j}{\mathrm{Im} \varphi_j(z)}\right)^{\gamma_j}\left|f(\varphi(z))\prod_{k=1}^n (\mathrm{Im} \varphi_k(z))^{\gamma_k}-f(\psi(z))\prod_{k=1}^n (\mathrm{Im} \psi_k(z))^{\gamma_k}\right|
		\end{align*}
	By lemma \ref{jiaquan hanshu mi cha}, we have
	\begin{align*}
	&h_1\prod_{j=1}^n \left(\frac{\mathrm{Im} z_j}{\mathrm{Im} \psi_j(z)}\right)^{\gamma_j}\left|f(\varphi(z))\prod_{k=1}^n (\mathrm{Im} \varphi_k(z))^{\gamma_k}-f(\psi(z))\prod_{k=1}^n (\mathrm{Im} \psi_k(z))^{\gamma_k}\right| \\
	&+h_2\prod_{j=1}^n \left(\frac{\mathrm{Im} z_j}{\mathrm{Im} \varphi_j(z)}\right)^{\gamma_j}\left|f(\varphi(z))\prod_{k=1}^n (\mathrm{Im} \varphi_k(z))^{\gamma_k}-f(\psi(z))\prod_{k=1}^n (\mathrm{Im} \psi_k(z))^{\gamma_k}\right| \\
	\lesssim&S_{\Omega_1\cup\Omega_2,f}^{\gamma}\rho(z)\left(h_2\prod_{k=1}^n \left(\frac{\mathrm{Im} z_k}{\mathrm{Im} \varphi_k(z)}\right)^{\gamma_k}+h_1\prod_{k=1}^n \left(\frac{\mathrm{Im} z_k}{\mathrm{Im} \psi_k(z)}\right)^{\gamma_k}\right).
	\end{align*}
	On the other hand, by lemma \ref{micha youxian}, we have
	\begin{align*}
		&\left|f(\varphi(z))\right|\prod_{j=1}^n (\mathrm{Im} \varphi_j(z))^{\gamma_j}\left|\prod_{k=1}^n \left(\frac{\mathrm{Im} z_k}{\mathrm{Im} \varphi_k(z)}\right)^{\gamma_k}-\prod_{k=1}^n \left(\frac{\mathrm{Im} z_k}{\mathrm{Im} \psi_k(z)}\right)^{\gamma_k}\right| \\
		\lesssim&S_{\Omega_1,f}^{\gamma}\rho(z)\left(\prod_{k=1}^n \left(\frac{\mathrm{Im} z_k}{\mathrm{Im} \varphi_k(z)}\right)^{\gamma_k}+\prod_{k=1}^n \left(\frac{\mathrm{Im} z_k}{\mathrm{Im} \psi_k(z)}\right)^{\gamma_k}\right).
	\end{align*}
	Similarly, we have
	\begin{align*}
		&\left|f(\psi(z))\right|\prod_{j=1}^n (\mathrm{Im} \psi_j(z))^{\gamma_j}\left|\prod_{k=1}^n \left(\frac{\mathrm{Im} z_k}{\mathrm{Im} \varphi_k(z)}\right)^{\gamma_k}-\prod_{k=1}^n \left(\frac{\mathrm{Im} z_k}{\mathrm{Im} \psi_k(z)}\right)^{\gamma_k}\right| \\
		\lesssim&S_{\Omega_1,f}^{\gamma}\rho(z)\left(\prod_{k=1}^n \left(\frac{\mathrm{Im} z_k}{\mathrm{Im} \varphi_k(z)}\right)^{\gamma_k}+\prod_{k=1}^n \left(\frac{\mathrm{Im} z_k}{\mathrm{Im} \psi_k(z)}\right)^{\gamma_k}\right).
	\end{align*}
	Combining these, we obtain
		\begin{align*}
		&\left|f(\varphi(z))-f(\psi(z))\right|\prod_{k=1}^n (\mathrm{Im} z_k)^{\gamma_k}\\
		\lesssim& \left(h_1 S_{\Omega_1,f}^{\gamma}+h_2 S_{\Omega_2,f}^{\gamma}\right)\rho(z)\left(\prod_{k=1}^n \left(\frac{\mathrm{Im} z_k}{\mathrm{Im} \varphi_k(z)}\right)^{\gamma_k}+\prod_{k=1}^n \left(\frac{\mathrm{Im} z_k}{\mathrm{Im} \psi_k(z)}\right)^{\gamma_k}\right) \\
		&+S_{\Omega_1\cup\Omega_2,f}^{\gamma}\rho(z)\left(h_2\prod_{k=1}^n \left(\frac{\mathrm{Im} z_k}{\mathrm{Im} \varphi_k(z)}\right)^{\gamma_k}+h_1\prod_{k=1}^n \left(\frac{\mathrm{Im} z_k}{\mathrm{Im} \psi_k(z)}\right)^{\gamma_k}\right).
	\end{align*}
\end{proof}

The following lemma is just a slight modification of the lemma $2.6$ in \cite{main}.
\begin{lemma}\label{compact}
	Let $\gamma\in\mathbb{R}_+^n$, and let $T$ be a complex linear combination of composition operators. Then $T: \mathcal{A}^{-\gamma}(\mathcal{H}^n) \to \mathcal{A}^{-\gamma}(\mathcal{H}^n)$ is compact if and only if $T: \mathcal{A}^{-\gamma}(\mathcal{H}^n) \to \mathcal{A}^{-\gamma}(\mathcal{H}^n)$ is bounded and for any bounded sequence in $\mathcal{A}^{-\gamma}(\mathcal{H}^n)$ that converges uniformly on compact subsets of $\mathcal{H}^n$ to zero, say a sequence $\{f^j\}$, we have $\|T f^j\|_{\mathcal{A}^{-\gamma}(\mathcal{H}^n)} \to 0$ as $j \to \infty$.
\end{lemma}

\section{Boundedness for difference of composition operators over $\mathcal{A}^{-\gamma}(\mathcal{H}^n)$}
\ \ \ \
In this section, we prove the first main theorem of this paper. We establish the following theorem, generalizing the conclusion of \cite{main}:

\begin{theorem}\label{thm1}
	Suppose $\varphi, \psi\in\mathcal{S}(\mathcal{H}^n)$, $\gamma\in\mathbb{R}_+^n$, then the following are equivalent:
\\
$(1)$ $\sup\limits_{z\in \mathcal{H}^n}\left(\prod\limits_{k=1}^n \left(\frac{\mathrm{Im} z_k}{\mathrm{Im} \varphi_k(z)}\right)^{\gamma_k}+\prod\limits_{k=1}^n \left(\frac{\mathrm{Im} z_k}{\mathrm{Im} \psi_k(z)}\right)^{\gamma_k}\right)\rho(z)<\infty$;
\\
$(2)$ $C_{\varphi}-C_{\psi}:\mathcal{A}^{-\gamma}(\mathcal{H}^n)\to\mathcal{A}^{-\gamma}(\mathcal{H}^n)$ is bounded.
\end{theorem}

\begin{proof}
	We first prove that $(1)\Rightarrow(2)$. Suppose that
	$$\sup\limits_{z\in \mathcal{H}^n}\left(\prod\limits_{k=1}^n \left(\frac{\mathrm{Im} z_k}{\mathrm{Im} \varphi_k(z)}\right)^{\gamma_k}+\prod\limits_{k=1}^n \left(\frac{\mathrm{Im} z_k}{\mathrm{Im} \psi_k(z)}\right)^{\gamma_k}\right)\rho(z)<\infty,$$
	Note that
	\begin{align*}
	&\|(C_{\varphi}-C_{\psi})f\|_{\mathcal{A}^{-\gamma}(\mathcal{H}^n)} \\ =&\sup_{z\in\mathcal{H}^n}\left|f(\varphi(z))-f(\psi(z))\right|\prod_{k=1}^n (\mathrm{Im} z_k)^{\gamma_k} \\
	=&\sup_{z\in\mathcal{H}^n}\left|f(\varphi(z))\prod_{k=1}^n (\mathrm{Im}\varphi_k(z))^{\gamma_k}\prod_{k=1}^n \bigg(\frac{\mathrm{Im} z_k}{\mathrm{Im}\varphi_k(z)}\bigg)^{\gamma_k}\right.\\
	&-\left.f(\psi(z))\prod_{k=1}^n (\mathrm{Im}\psi_k(z))^{\gamma_k}\prod_{k=1}^n \bigg(\frac{\mathrm{Im} z_k}{\mathrm{Im}\psi_k(z)}\bigg)^{\gamma_k}\right| \\
	\leqslant&\sup_{z\in\mathcal{H}^n}|f(\varphi(z))|\prod_{k=1}^n (\mathrm{Im}\varphi_k(z))^{\gamma_k}\left|\prod_{k=1}^n \left(\frac{\mathrm{Im} z_k}{\mathrm{Im}\varphi_k(z)}\right)^{\gamma_k}-\prod_{k=1}^n \left(\frac{\mathrm{Im} z_k}{\mathrm{Im}\psi_k(z)}\right)^{\gamma_k}\right| \\
	&+\sup_{z\in\mathcal{H}^n}\left|f(\varphi(z))\prod_{k=1}^n (\mathrm{Im}\varphi_k(z))^{\gamma_k}-f(\psi(z))\prod_{k=1}^n (\mathrm{Im}\psi_k(z))^{\gamma_k} \right|\prod_{k=1}^n \left(\frac{\mathrm{Im} z_k}{\mathrm{Im}\psi_k(z)}\right)^{\gamma_k}.
	\end{align*}
	From lemma \ref{jiaquan hanshu mi cha}, it is clear that
	\begin{align*}
		&\sup_{z\in\mathcal{H}^n}\left|f(\varphi(z))\prod_{k=1}^n (\mathrm{Im}\varphi_k(z))^{\gamma_k}-f(\psi(z))\prod_{k=1}^n (\mathrm{Im}\psi_k(z))^{\gamma_k} \right| \\
		\lesssim&\|f\|_{\mathcal{A}^{-\gamma}(\mathcal{H}^n)}\sup_{z\in\mathcal{H}^n}\rho(z)\prod_{k=1}^n \left(\frac{\mathrm{Im} z_k}{\mathrm{Im}\psi_k(z)}\right)^{\gamma_k}.
	\end{align*}
	Then by lemma \ref{micha youxian}, we deduce that
\begin{align*}
	&\|(C_{\varphi}-C_{\psi})f\|_{\mathcal{A}^{-\gamma}(\mathcal{H}^n)} \\
	\lesssim&\|f\|_{\mathcal{A}^{-\gamma}(\mathcal{H}^n)}\sup_{z\in\mathcal{H}^n}\left|\prod_{k=1}^n \left(\frac{\mathrm{Im} z_k}{\mathrm{Im}\varphi_k(z)}\right)^{\gamma_k}-\prod_{k=1}^n \left(\frac{\mathrm{Im} z_k}{\mathrm{Im}\psi_k(z)}\right)^{\gamma_k}\right|\\
	&+\|f\|_{\mathcal{A}^{-\gamma}(\mathcal{H}^n)}\sup_{z\in\mathcal{H}^n}\rho(z)\prod_{k=1}^n \left(\frac{\mathrm{Im} z_k}{\mathrm{Im}\psi_k(z)}\right)^{\gamma_k} \\
	\lesssim&\|f\|_{\mathcal{A}^{-\gamma}(\mathcal{H}^n)}\sup_{z\in\mathcal{H}^n}\rho(z)\left(\prod_{k=1}^n \left(\frac{\mathrm{Im} z_k}{\mathrm{Im}\varphi_k(z)}\right)^{\gamma_k}+\prod_{k=1}^n \left(\frac{\mathrm{Im} z_k}{\mathrm{Im}\psi_k(z)}\right)^{\gamma_k}\right) \\
	<&\infty.
\end{align*}

For the proof of $(2)\Rightarrow(1)$, assume that $C_{\varphi}-C_{\psi}$ is a bounded operator. Notice that for any $m\in\{1,\cdots,n\}$, if we consider the function $g_{a,m}$ as defined in (\ref{g_a,m}), we have
\begin{align*}
	&\|(C_{\varphi}-C_{\psi})g_{\varphi(z),m}\|_{\mathcal{A}^{-\gamma}(\mathcal{H}^n)} \\
	=&\sup_{w\in\mathcal{H}^n}|g_{\varphi(z),m}(\varphi(w))-g_{\varphi(z),m}(\psi(w))|\prod_{k=1}^{n}(\mathrm{Im} w_k)^{\gamma_k} \\
	\geqslant&|g_{\varphi(z),m}(\varphi(z))-g_{\varphi(z),m}(\psi(z))|\prod_{k=1}^{n}(\mathrm{Im} z_k)^{\gamma_k} \\
	=&4^{|\gamma|}\left|\frac{\psi_m(z)-\varphi_m(z)}{\psi_m(z)-\overline{\varphi_m(z)}}\cdot\prod_{k=1}^n \frac{(\mathrm{Im} \varphi_k(z))^{\gamma_k}}{(\psi_k(z)-\overline{\varphi_k(z)})^{2\gamma_k}}\right|\prod_{k=1}^{n}(\mathrm{Im} z_k)^{\gamma_k} \\
	=&4^{|\gamma|}\rho_m(z)\prod_{k=1}^n \left(\frac{\mathrm{Im} z_k}{\mathrm{Im} \varphi_k(z)}\right)^{\gamma_k}\left|\frac{\mathrm{Im} \varphi_k(z)}{\psi_k(z)-\overline{\varphi_k(z)}}\right|^{2\gamma_k}.
\end{align*}
Summing up over $m$ from $1$ to $n$, we obtain
\begin{align}\label{g bound}
	&\sum_{m=1}^{n}\|(C_{\varphi}-C_{\psi})g_{\varphi(z),m}\|_{\mathcal{A}^{-\gamma}(\mathcal{H}^n)}\notag \\
	=&4^{|\gamma|}\left(\sum_{m=1}^{n}\rho_m(z)\right)\prod_{k=1}^n \left(\frac{\mathrm{Im} z_k}{\mathrm{Im} \varphi_k(z)}\right)^{\gamma_k}\left|\frac{\mathrm{Im} \varphi_k(z)}{\psi_k(z)-\overline{\varphi_k(z)}}\right|^{2\gamma_k} \notag\\
	\geqslant&4^{|\gamma|}\max_{1\leqslant m\leqslant n}\rho_m(z)\prod_{k=1}^n \left(\frac{\mathrm{Im} z_k}{\mathrm{Im} \varphi_k(z)}\right)^{\gamma_k}\left|\frac{\mathrm{Im} \varphi_k(z)}{\psi_k(z)-\overline{\varphi_k(z)}}\right|^{2\gamma_k} \notag\\	
	=&4^{|\gamma|}\rho(z)\prod_{k=1}^n \left(\frac{\mathrm{Im} z_k}{\mathrm{Im} \varphi_k(z)}\right)^{\gamma_k}\left|\frac{\mathrm{Im} \varphi_k(z)}{\psi_k(z)-\overline{\varphi_k(z)}}\right|^{2\gamma_k}.		
\end{align}
	On the other hand, by the definition in (\ref{f_a}) of the function $f_a$, we have
\begin{align}\label{f bound 1}
	&\prod_{k=1}^n \left(\frac{\mathrm{Im} z_k}{\mathrm{Im}\varphi_k(z)}\right)^{\gamma_k}-4^{|\gamma|}\prod_{k=1}^n \left(\frac{\mathrm{Im} z_k}{\mathrm{Im}\varphi_k(z)}\cdot\left|\frac{\mathrm{Im} \varphi_k(z)}{\psi_k(z)-\overline{\varphi_k(z)}}\right|^2\right)^{\gamma_k}\notag \\
	\leqslant&\left|\prod_{k=1}^n \left(\frac{(2i)^2 \mathrm{Im}\varphi_k(z)}{(\varphi_k(z)-\overline{\varphi_k(z)})^2}\right)^{\gamma_k}-\prod_{k=1}^n \left(\frac{(2i)^2 \mathrm{Im} \varphi_k(z)}{(\psi_k(z)-\overline{\varphi_k(z)})^2}\right)^{\gamma_k}\right|\prod_{k=1}^{n}(\mathrm{Im} z_k)^{\gamma_k}\notag \\
	\leqslant&\sup_{w\in\mathcal{H}^n}\left|\prod_{k=1}^n \left(\frac{(2i)^2 \mathrm{Im}\varphi_k(z)}{(\varphi_k(w)-\overline{\varphi_k(z)})^2}\right)^{\gamma_k}-\prod_{k=1}^n \left(\frac{(2i)^2 \mathrm{Im} \varphi_k(z)}{(\psi_k(w)-\overline{\varphi_k(z)})^2}\right)^{\gamma_k}\right|\prod_{k=1}^{n}(\mathrm{Im} w_k)^{\gamma_k}\notag \\
	\leqslant&\|(C_\varphi-C_\psi)f_{\varphi(z)}\|_{\mathcal{A}^{-\gamma}(\mathcal{H}^n)}.
\end{align}
	For the same reason, 
	\begin{align}\label{f bound 2}
		&\prod_{k=1}^n \left(\frac{\mathrm{Im} z_k}{\mathrm{Im}\psi_k(z)}\right)^{\gamma_k}-4^{|\gamma|}\prod_{k=1}^n \left(\frac{\mathrm{Im} z_k}{\mathrm{Im}\psi_k(z)}\cdot\left|\frac{\mathrm{Im} \psi_k(z)}{\varphi_k(z)-\overline{\psi_k(z)}}\right|^2\right)^{\gamma_k}\notag \\
		\leqslant&\|(C_\varphi-C_\psi)f_{\psi(z)}\|_{\mathcal{A}^{-\gamma}(\mathcal{H}^n)}.
	\end{align}
	Therefore, 
\begin{align*}
	&\rho(z)\left(\prod_{k=1}^n \left(\frac{\mathrm{Im} z_k}{\mathrm{Im}\varphi_k(z)}\right)^{\gamma_k}+\prod_{k=1}^n \left(\frac{\mathrm{Im} z_k}{\mathrm{Im}\psi_k(z)}\right)^{\gamma_k}\right) \\
	\leqslant&\rho(z)\prod_{k=1}^n \left(\frac{\mathrm{Im} z_k}{\mathrm{Im}\varphi_k(z)}\right)^{\gamma_k}-4^{|\gamma|}\rho(z)\prod_{k=1}^n \left(\frac{\mathrm{Im} z_k}{\mathrm{Im}\varphi_k(z)}\cdot\left|\frac{\mathrm{Im} \varphi_k(z)}{\psi_k(z)-\overline{\varphi_k(z)}}\right|^2\right)^{\gamma_k} \\
	&+\rho(z)\prod_{k=1}^n \left(\frac{\mathrm{Im} z_k}{\mathrm{Im}\psi_k(z)}\right)^{\gamma_k}-4^{|\gamma|}\rho(z)\prod_{k=1}^n \left(\frac{\mathrm{Im} z_k}{\mathrm{Im}\psi_k(z)}\cdot\left|\frac{\mathrm{Im} \psi_k(z)}{\varphi_k(z)-\overline{\psi_k(z)}}\right|^2\right)^{\gamma_k} \\
	&+4^{|\gamma|}\rho(z)\prod_{k=1}^n \left(\frac{\mathrm{Im} z_k}{\mathrm{Im} \varphi_k(z)}\right)^{\gamma_k}\left|\frac{\mathrm{Im} \varphi_k(z)}{\psi_k(z)-\overline{\varphi_k(z)}}\right|^{2\gamma_k} \\
	&+4^{|\gamma|}\rho(z)\prod_{k=1}^n \left(\frac{\mathrm{Im} z_k}{\mathrm{Im} \psi_k(z)}\right)^{\gamma_k}\left|\frac{\mathrm{Im} \psi_k(z)}{\varphi_k(z)-\overline{\psi_k(z)}}\right|^{2\gamma_k} \\
	\leqslant&\rho(z)\|(C_\varphi-C_\psi)f_{\varphi(z)}\|_{\mathcal{A}^{-\gamma}(\mathcal{H}^n)}+\rho(z)\|(C_\varphi-C_\psi)f_{\psi(z)}\|_{\mathcal{A}^{-\gamma}(\mathcal{H}^n)} \\
	&+\sum_{m=1}^{n}\|(C_{\varphi}-C_{\psi})g_{\varphi(z),m}\|_{\mathcal{A}^{-\gamma}(\mathcal{H}^n)}+\sum_{m=1}^{n}\|(C_{\varphi}-C_{\psi})g_{\psi(z),m}\|_{\mathcal{A}^{-\gamma}(\mathcal{H}^n)},
\end{align*}
Thus, from the boundedness of $C_{\varphi}-C_{\psi}$,
\begin{align*}
	&\sup\limits_{z\in \mathcal{H}^n}\rho(z)\left(\prod_{k=1}^n \left(\frac{\mathrm{Im} z_k}{\mathrm{Im}\varphi_k(z)}\right)^{\gamma_k}+\prod_{k=1}^n \left(\frac{\mathrm{Im} z_k}{\mathrm{Im}\psi_k(z)}\right)^{\gamma_k}\right) \\
	\leqslant&\sup_{z\in \mathcal{H}^n}\rho(z)\|(C_\varphi-C_\psi)f_{\varphi(z)}\|_{\mathcal{A}^{-\gamma}(\mathcal{H}^n)}+\sup_{z\in \mathcal{H}^n}\rho(z)\|(C_\varphi-C_\psi)f_{\psi(z)}\|_{\mathcal{A}^{-\gamma}(\mathcal{H}^n)} \\
	&+\sum_{m=1}^{n}\sup_{z\in \mathcal{H}^n}\|(C_{\varphi}-C_{\psi})g_{\varphi(z),m}\|_{\mathcal{A}^{-\gamma}(\mathcal{H}^n)}+\sum_{m=1}^{n}\sup_{z\in \mathcal{H}^n}\|(C_{\varphi}-C_{\psi})g_{\psi(z),m}\|_{\mathcal{A}^{-\gamma}(\mathcal{H}^n)} \\
	\leqslant&\sup_{a\in \mathcal{H}^n}\rho(z)\|(C_\varphi-C_\psi)f_a\|_{\mathcal{A}^{-\gamma}(\mathcal{H}^n)}+\sup_{a\in \mathcal{H}^n}\rho(z)\|(C_\varphi-C_\psi)f_a\|_{\mathcal{A}^{-\gamma}(\mathcal{H}^n)} \\
	&+\sum_{m=1}^{n}\sup_{a\in \mathcal{H}^n}\|(C_{\varphi}-C_{\psi})g_{a,m}\|_{\mathcal{A}^{-\gamma}(\mathcal{H}^n)}+\sum_{m=1}^{n}\sup_{a\in \mathcal{H}^n}\|(C_{\varphi}-C_{\psi})g_{a,m}\|_{\mathcal{A}^{-\gamma}(\mathcal{H}^n)} \\
	\lesssim&\rho(z)\sup_{a\in \mathcal{H}^n}\|f_a\|_{\mathcal{A}^{-\gamma}(\mathcal{H}^n)}+\sum_{m=1}^{n}\sup_{a\in \mathcal{H}^n}\|g_{a,m}\|_{\mathcal{A}^{-\gamma}(\mathcal{H}^n)}<\infty,
\end{align*}
which completes the proof.
\end{proof}

Furthermore, we find that the first part of theorem \ref{thm1} also implies that $C_{\varphi}-C_{\psi}$ is an absolutely summing operator and hence we have the following corollary. 
\begin{corollary}
	Suppose $\varphi, \psi\in\mathcal{S}(\mathcal{H}^n)$, $\gamma\in\mathbb{R}_+^n$, $1\leq p<\infty$, the operator $C_{\varphi}-C_{\psi}:\mathcal{A}^{-\gamma}(\mathcal{H}^n)\to\mathcal{A}^{-\gamma}(\mathcal{H}^n)$, then the following are equivalent:
	\\
	$(1)$ $C_{\varphi}-C_{\psi}$ is an absolutely summing operator;
	\\
	$(2)$ $C_{\varphi}-C_{\psi}$ is a $p$-summing operator;
	\\
	$(3)$ $C_{\varphi}-C_{\psi}$ is bounded.
\end{corollary}

\begin{proof}
	As we have shown in section 2, $(1)\Rightarrow(2)$ and $(2)\Rightarrow(3)$ are well-known results of $p$-summing operators for Banach space, thus it remains to prove that $(3)\Rightarrow(1)$ holds.
	
	For any $\xi\in\mathcal{H}^n$, we consider a bounded linear functional $u_{\xi}\in\mathcal{A}^{-\gamma}(\mathcal{H}^n)^*$ given by $u_{\xi}(f)=f(\xi)\prod\limits_{k=1}^{n}(\mathrm{Im} \xi_k)^{\gamma_k}$
	
	Since
	$$\left|u_{\xi}(f)\right|=\left|f(\xi)\prod\limits_{k=1}^{n}(\mathrm{Im} \xi_k)^{\gamma_k}\right|\leqslant\|f\|_{\mathcal{A}^{-\gamma}(\mathcal{H}^n)},$$
	we have $\|u_{\xi}\|\leqslant 1$. On the other hand, 
	$$\left\|f_{\xi}\right\|_{\mathcal{A}^{-\gamma}(\mathcal{H}^n)}=\sup_{z\in\mathcal{H}^n}|f_{\xi}(z)|\prod_{k=1}^{n}(\mathrm{Im} z_k)^{\gamma_k}=|f_{\xi}(\xi)|\prod_{k=1}^{n}(\mathrm{Im} \xi_k)^{\gamma_k}=\left|u_{\xi}(f_{\xi})\right|.$$
	
	It follows that $\|u_{\xi}\|=1$ and $u_{\xi}\in B_{\mathcal{A}^{-\gamma}(\mathcal{H}^n)^*}$, so if we suppose that $f^1,f^2,\cdots,$ $f^m\in\mathcal{A}^{-\gamma}(\mathcal{H}^n)$ then
\begin{align*}
	\sup_{u\in B_{\mathcal{A}^{-\gamma}(\mathcal{H}^n)^*}}\sum_{j\leqslant m}|u(f^j)|\geqslant&\sup_{\xi\in\mathcal{H}^n} \sum_{j \leq m} |u_{\xi}(f^j)| \\
	=&\sup_{\xi\in\mathcal{H}^n} \sum_{j \leq m} \left|f^j(\xi)\prod_{k=1}^{n}(\mathrm{Im} \xi_k)^{\gamma_k}\right| \\
	=&\sum_{j \leq m}\left(\sup_{\xi\in\mathcal{H}^n}  \left|f^j(\xi)\right|\prod_{k=1}^{n}(\mathrm{Im} \xi_k)^{\gamma_k}\right) \\
	=&\sum_{j \leq m}\|f^j\|_{\mathcal{A}^{-\gamma}(\mathcal{H}^n)}.
\end{align*}
	With this and the calculation appeared in the proof of $(1)\Rightarrow(2)$ of theorem \ref{thm1}, we obtain
\begin{align*}
	&\sum_{j \leq m} \|(C_{\varphi}-C_{\psi})f^j \|_{\mathcal{A}^{-\gamma}(\mathcal{H}^n)} \\
	\lesssim&\sum_{j \leq m}\left(\|f^j\|_{\mathcal{A}^{-\gamma}(\mathcal{H}^n)}\sup_{z\in\mathcal{H}^n}\rho(z)\left(\prod_{k=1}^n \left(\frac{\mathrm{Im} z_k}{\mathrm{Im}\varphi_k(z)}\right)^{\gamma_k}+\prod_{k=1}^n \left(\frac{\mathrm{Im} z_k}{\mathrm{Im}\psi_k(z)}\right)^{\gamma_k}\right)\right) \\
	=&\sup_{z\in\mathcal{H}^n}\rho(z)\left(\prod_{k=1}^n \left(\frac{\mathrm{Im} z_k}{\mathrm{Im}\varphi_k(z)}\right)^{\gamma_k}+\prod_{k=1}^n \left(\frac{\mathrm{Im} z_k}{\mathrm{Im}\psi_k(z)}\right)^{\gamma_k}\right)\sum_{j \leq m}\|f^j\|_{\mathcal{A}^{-\gamma}(\mathcal{H}^n)} \\
	\leqslant&\sup_{z\in\mathcal{H}^n}\rho(z)\left(\prod_{k=1}^n \left(\frac{\mathrm{Im} z_k}{\mathrm{Im}\varphi_k(z)}\right)^{\gamma_k}+\prod_{k=1}^n \left(\frac{\mathrm{Im} z_k}{\mathrm{Im}\psi_k(z)}\right)^{\gamma_k}\right)\sup_{u\in B_{\mathcal{A}^{-\gamma}(\mathcal{H}^n)^*}}\sum_{j\leqslant m}|u(f^j)|.
\end{align*}
	
	From Theorem \ref{thm1} we knew that if $C_{\varphi}-C_{\psi}$ is bounded then 
	$$\sup_{z\in\mathcal{H}^n}\rho(z)\left(\prod_{k=1}^n \left(\frac{\mathrm{Im} z_k}{\mathrm{Im}\varphi_k(z)}\right)^{\gamma_k}+\prod_{k=1}^n \left(\frac{\mathrm{Im} z_k}{\mathrm{Im}\psi_k(z)}\right)^{\gamma_k}\right)<\infty,$$
	which means that
	$$\sum_{j \leq m} \|(C_{\varphi}-C_{\psi})f^j \|_{\mathcal{A}^{-\gamma}(\mathcal{H}^n)}\lesssim\sup_{u\in B_{\mathcal{A}^{-\gamma}(\mathcal{H}^n)^*}}\sum_{j\leqslant m}|u(f^j)|,$$
	and thus we have shown that $C_{\varphi}-C_{\psi}$ is an absolutely summing operator.

\end{proof}

\section{Compactness for difference of composition operators over $\mathcal{A}^{-\gamma}(\mathcal{H}^n)$}
\ \ \ \
Next, we will prove in this section the second main theorem of this paper, that is, the characterization of the compactness of $C_{\varphi}-C_{\psi}$.

\begin{theorem}
	Suppose $\varphi, \psi\in\mathcal{S}(\mathcal{H}^n)$, $\gamma\in\mathbb{R}_+^n$, the operator $C_{\varphi}-C_{\psi}:\mathcal{A}^{-\gamma}(\mathcal{H}^n)\to\mathcal{A}^{-\gamma}(\mathcal{H}^n)$, then the following are equivalent:
	\\
	$(1)$ $C_{\varphi}-C_{\psi}$ is compact;
	\\
	$(2)$ $\lim\limits_{\varphi(z)\to\partial\widehat{\mathcal{H}^n}}\rho(z)\prod\limits_{k=1}^n \left(\frac{\mathrm{Im} z_k}{\mathrm{Im} \varphi_k(z)}\right)^{\gamma_k}+\lim\limits_{\psi(z)\to\partial\widehat{\mathcal{H}^n}}\rho(z)\prod\limits_{k=1}^n \left(\frac{\mathrm{Im} z_k}{\mathrm{Im} \psi_k(z)}\right)^{\gamma_k}=0$.
\end{theorem}

\begin{proof}
	Our discussion begin with the proof of $(1)\Rightarrow(2)$. To reach the result we only need to show that if $C_{\varphi}-C_{\psi}$ is compact then
	$$\lim\limits_{\varphi(z)\to\partial\widehat{\mathcal{H}^n}}\rho(z)\prod\limits_{k=1}^n \left(\frac{\mathrm{Im} z_k}{\mathrm{Im} \varphi_k(z)}\right)^{\gamma_k}=0\quad\text{and}\quad\lim\limits_{\psi(z)\to\partial\widehat{\mathcal{H}^n}}\rho(z)\prod\limits_{k=1}^n \left(\frac{\mathrm{Im} z_k}{\mathrm{Im} \psi_k(z)}\right)^{\gamma_k}=0.$$
	
	We assume that $C_{\varphi}-C_{\psi}:\mathcal{A}^{-\gamma}(\mathcal{H}^n)\to\mathcal{A}^{-\gamma}(\mathcal{H}^n)$ is a compact operator, then by lemma \ref{compact} for any bounded sequence $\{f^j\}$ in $\mathcal{A}^{-\gamma}(\mathcal{H}^n)$ that converges uniformly to $0$ on each compact subset of $\mathcal{H}^n$, we have
	$$\limsup_{j\to\infty}\|(C_{\varphi}-C_{\psi})f^j \|_{\mathcal{A}^{-\gamma}(\mathcal{H}^n)}=0.$$
	In particularly, for the test functions $f_a$ and $g_{a,m}$ given in (\ref{f_a}) and (\ref{g_a,m}), 
	\begin{align}
		\limsup_{a\to\partial\widehat{\mathcal{H}^n}}\|(C_{\varphi}-C_{\psi})f_a \|_{\mathcal{A}^{-\gamma}(\mathcal{H}^n)}=0, \\
		\limsup_{a\to\partial\widehat{\mathcal{H}^n}}\|(C_{\varphi}-C_{\psi})g_{a,m} \|_{\mathcal{A}^{-\gamma}(\mathcal{H}^n)}=0.
	\end{align}

	Using (\ref{g bound}) and (\ref{f bound 1}) we obtain
	\begin{align*}
		&\rho(z)\prod_{k=1}^n \left(\frac{\mathrm{Im} z_k}{\mathrm{Im}\varphi_k(z)}\right)^{\gamma_k} \\
		\leqslant&\rho(z)\prod_{k=1}^n \left(\frac{\mathrm{Im} z_k}{\mathrm{Im}\varphi_k(z)}\right)^{\gamma_k}-4^{|\gamma|}\rho(z)\prod_{k=1}^n \left(\frac{\mathrm{Im} z_k}{\mathrm{Im}\varphi_k(z)}\cdot\left|\frac{\mathrm{Im} \varphi_k(z)}{\psi_k(z)-\overline{\varphi_k(z)}}\right|^2\right)^{\gamma_k} \\
		&+4^{|\gamma|}\rho(z)\prod_{k=1}^n \left(\left(\frac{\mathrm{Im} z_k}{\mathrm{Im} \varphi_k(z)}\right)^{\gamma_k}\left|\frac{\mathrm{Im} \varphi_k(z)}{\psi_k(z)-\overline{\varphi_k(z)}}\right|^{2\gamma_k}\right) \\ 
		\leqslant&\rho(z)\|(C_\varphi-C_\psi)f_{\varphi(z)}\|_{\mathcal{A}^{-\gamma}(\mathcal{H}^n)}+\sum_{m=1}^{n}\|(C_{\varphi}-C_{\psi})g_{\varphi(z),m}\|_{\mathcal{A}^{-\gamma}(\mathcal{H}^n)},
	\end{align*}
	
	Hence,
	\begin{align*}
		&\lim_{\varphi(z)\to\partial\widehat{\mathcal{H}^n}}\rho(z)\prod_{k=1}^n \left(\frac{\mathrm{Im} z_k}{\mathrm{Im}\varphi_k(z)}\right)^{\gamma_k} \\
		\lesssim&\limsup_{\varphi(z)\to\partial\widehat{\mathcal{H}^n}}\rho(z)\|(C_\varphi-C_\psi)f_{\varphi(z)}\|_{\mathcal{A}^{-\gamma}(\mathcal{H}^n)}\\
		&+\sum_{m=1}^{n}\limsup_{\varphi(z)\to\partial\widehat{\mathcal{H}^n}}\|(C_{\varphi}-C_{\psi})g_{\varphi(z),m}\|_{\mathcal{A}^{-\gamma}(\mathcal{H}^n)} \\
		\lesssim&\limsup_{a\to\partial\widehat{\mathcal{H}^n}}\rho(z)\|(C_\varphi-C_\psi)f_a\|_{\mathcal{A}^{-\gamma}(\mathcal{H}^n)}\\
		&+\sum_{m=1}^{n}\limsup_{a\to\partial\widehat{\mathcal{H}^n}}\|(C_{\varphi}-C_{\psi})g_{a,m}\|_{\mathcal{A}^{-\gamma}(\mathcal{H}^n)} \\
		=&0.
	\end{align*}
Similarly, we have
$$\lim_{\psi(z)\to\partial\widehat{\mathcal{H}^n}}\rho(z)\prod_{k=1}^n \left(\frac{\mathrm{Im} z_k}{\mathrm{Im}\psi_k(z)}\right)^{\gamma_k}=0.$$
By combining the two limits, we obtain
$$\lim_{\varphi(z)\to\partial\widehat{\mathcal{H}^n}}\rho(z)\prod\limits_{k=1}^n \left(\frac{\mathrm{Im} z_k}{\mathrm{Im} \varphi_k(z)}\right)^{\gamma_k}+\lim_{\psi(z)\to\partial\widehat{\mathcal{H}^n}}\rho(z)\prod\limits_{k=1}^n \left(\frac{\mathrm{Im} z_k}{\mathrm{Im} \psi_k(z)}\right)^{\gamma_k}=0.$$
	
	To prove $(2)\Rightarrow(1)$, we assume that $(2)$ holds. From known results, $C_{\varphi}-C_{\psi}$ is bounded on $A^{-\gamma}(\mathcal{H}^n)$, thus
	$$R_{\varphi,\psi}:=\sup_{z\in\mathcal{H}^n}\rho(z)\left(\prod_{k=1}^n \left(\frac{\mathrm{Im} z_k}{\mathrm{Im}\varphi_k(z)}\right)^{\gamma_k}+\prod_{k=1}^n \left(\frac{\mathrm{Im} z_k}{\mathrm{Im}\psi_k(z)}\right)^{\gamma_k}\right)<\infty.$$ 
	Let $\{f^j\} \subset A^{-\gamma}(\mathcal{H}^n)$ be any sequence such that $\| f^j \|_{A^{-\gamma}(\mathcal{H}^n)} \leq 1$ and $\{f^j\}$ converges uniformly to $0$ on any compact subset of $\mathcal{H}^n$ as $j\to\infty$. Moreover, for any $\varepsilon>0$ there is a $J>0$ such that
	$$S_{D,f^j}^{\gamma}<\frac{\varepsilon}{2R_{\varphi,\psi}}$$
	holds on any compact subset $D\subset\mathcal{H}^n$ whenever $j>J$.
	
	To prove that $C_{\varphi}-C_{\psi}$ is compact, by lemma \ref{compact} it suffices to prove that $\| (C_{\varphi}-C_{\psi)}f^j \|_{A^{-\gamma}(\mathcal{H}^n)} \to 0$ as $j \to \infty $. By conditions $(2)$, for any $\varepsilon > 0$ and $z\in\mathcal{H}^n$, there exists a compact set $\overline{\Omega} \subset\mathcal{H}^n$ where $\Omega$ is a bounded open subset of $\mathcal{H}^n$ such that $\varphi(z),\psi(z)\notin\partial\Omega$ and

\begin{equation}
	\rho(z)\prod_{k=1}^n \left(\frac{\mathrm{Im} z_k}{\mathrm{Im}\varphi_k(z)}\right)^{\gamma_k}<\frac{\varepsilon}{2},\quad \text{if}\quad\varphi(z) \in \mathcal{H}^n\setminus\overline{\Omega},
\end{equation}
\begin{equation}
	\rho(z)\prod_{k=1}^n \left(\frac{\mathrm{Im} z_k}{\mathrm{Im}\psi_k(z)}\right)^{\gamma_k}<\frac{\varepsilon}{2},\quad \text{if}\quad\psi(z) \in\mathcal{H}^n \setminus\overline{\Omega}.
\end{equation}

For arbitrary $\varepsilon>0$, using lemma \ref{limit}, we have
\begin{align}\label{saikou}
	&\left|f^j(\varphi(z))-f^j(\psi(z))\right|\prod_{k=1}^n (\mathrm{Im} z_k)^{\gamma_k} \notag\\
	\leqslant& C_{\gamma,\Omega}\left(\left(h_1 S_{\Omega_1,f^j}^{\gamma}+h_2 S_{\Omega_2,f^j}^{\gamma}\right)\rho(z)\left(\prod_{k=1}^n \left(\frac{\mathrm{Im} z_k}{\mathrm{Im} \varphi_k(z)}\right)^{\gamma_k}+\prod_{k=1}^n \left(\frac{\mathrm{Im} z_k}{\mathrm{Im} \psi_k(z)}\right)^{\gamma_k}\right)\right. \notag\\
	&+\left.S_{\Omega_1\cup\Omega_2,f^j}^{\gamma}\rho(z)\left(h_2\prod_{k=1}^n \left(\frac{\mathrm{Im} z_k}{\mathrm{Im} \varphi_k(z)}\right)^{\gamma_k}+h_1\prod_{k=1}^n \left(\frac{\mathrm{Im} z_k}{\mathrm{Im} \psi_k(z)}\right)^{\gamma_k}\right)\right) \notag\\
	\leqslant& C_{\gamma,\Omega}\left(\left(h_1 S_{\overline{\Omega_1},f^j}^{\gamma}+h_2 S_{\overline{\Omega_2},f^j}^{\gamma}\right)\rho(z)\left(\prod_{k=1}^n \left(\frac{\mathrm{Im} z_k}{\mathrm{Im} \varphi_k(z)}\right)^{\gamma_k}+\prod_{k=1}^n \left(\frac{\mathrm{Im} z_k}{\mathrm{Im} \psi_k(z)}\right)^{\gamma_k}\right)\right. \notag\\
	&+\left.S_{\overline{\Omega_1\cup\Omega_2},f^j}^{\gamma}\rho(z)\left(h_2\prod_{k=1}^n \left(\frac{\mathrm{Im} z_k}{\mathrm{Im} \varphi_k(z)}\right)^{\gamma_k}+h_1\prod_{k=1}^n \left(\frac{\mathrm{Im} z_k}{\mathrm{Im} \psi_k(z)}\right)^{\gamma_k}\right)\right)
\end{align}
where $\Omega_1,\Omega_2\in\{\Omega,\quad\mathcal{H}^n\setminus\overline{\Omega}\}$. The constant $C_{\gamma,\Omega}$ here is independent of the choice of $h_1$ and $h_2$. To finish the proof, we divide the remaining discussion into four cases.

\noindent $\bf{(a)}$ If $\varphi(z)\in\Omega,\psi(z)\in\Omega$, then by (\ref{saikou})
\begin{align*}
	&\left|f^j(\varphi(z))-f^j(\psi(z))\right|\prod_{k=1}^n (\mathrm{Im} z_k)^{\gamma_k} \\
	\leqslant& C_{\gamma,\Omega}\left(\left(h_1 S_{\overline{\Omega},f^j}^{\gamma}+h_2 S_{\overline{\Omega},f^j}^{\gamma}\right)\rho(z)\left(\prod_{k=1}^n \left(\frac{\mathrm{Im} z_k}{\mathrm{Im} \varphi_k(z)}\right)^{\gamma_k}+\prod_{k=1}^n \left(\frac{\mathrm{Im} z_k}{\mathrm{Im} \psi_k(z)}\right)^{\gamma_k}\right)\right. \\
	&+\left.S_{\overline{\Omega\cup\Omega},f^j}^{\gamma}\rho(z)\left(h_2\prod_{k=1}^n \left(\frac{\mathrm{Im} z_k}{\mathrm{Im} \varphi_k(z)}\right)^{\gamma_k}+h_1\prod_{k=1}^n \left(\frac{\mathrm{Im} z_k}{\mathrm{Im} \psi_k(z)}\right)^{\gamma_k}\right)\right) \\
	\leqslant& 2C_{\gamma,\Omega}S_{\overline{\Omega},f^j}^{\gamma}R_{\varphi,\psi}
\end{align*}
Since $\{f^j\}$ is uniformly convergent to $0$ on $\overline{\Omega}$ as $j\to\infty$, as mentioned earlier in this paper, we can choose a $J>0$ such that $$S_{\overline{\Omega},f^j}^{\gamma}<\frac{\varepsilon}{2R_{\varphi,\psi}}$$
whenever $j>J$ and therefore
$$\left|f^j(\varphi(z))-f^j(\psi(z))\right|\prod_{k=1}^n (\mathrm{Im} z_k)^{\gamma_k}<C_{\gamma,\Omega}\varepsilon.$$

\noindent $\bf{(b)}$ If $\varphi(z)\in\Omega,\psi(z)\in\mathcal{H}^n \setminus\overline{\Omega}$,
\begin{align*}
	&\left|f^j(\varphi(z))-f^j(\psi(z))\right|\prod_{k=1}^n (\mathrm{Im} z_k)^{\gamma_k} \\
	\leqslant& C_{\gamma,\Omega}\left(\left(h_1 S_{\overline{\Omega},f^j}^{\gamma}+h_2 S_{\mathcal{H}^n \setminus\Omega,f^j}^{\gamma}\right)\rho(z)\left(\prod_{k=1}^n \left(\frac{\mathrm{Im} z_k}{\mathrm{Im} \varphi_k(z)}\right)^{\gamma_k}+\prod_{k=1}^n \left(\frac{\mathrm{Im} z_k}{\mathrm{Im} \psi_k(z)}\right)^{\gamma_k}\right)\right. \\
	&+\left.S_{\mathcal{H}^n,f^j}^{\gamma}\rho(z)\left(h_2\prod_{k=1}^n \left(\frac{\mathrm{Im} z_k}{\mathrm{Im} \varphi_k(z)}\right)^{\gamma_k}+h_1\prod_{k=1}^n \left(\frac{\mathrm{Im} z_k}{\mathrm{Im} \psi_k(z)}\right)^{\gamma_k}\right)\right) \\
	\leqslant&C_{\gamma,\Omega}\left(\left(h_1 S_{\overline{\Omega},f^j}^{\gamma}+h_2 S_{\mathcal{H}^n \setminus\Omega,f^j}^{\gamma}\right)R_{\varphi,\psi}\right. \\
	&+\left.\| f^j \|_{A^{-\gamma}(\mathcal{H}^n)}\rho(z)\left(h_2\prod_{k=1}^n \left(\frac{\mathrm{Im} z_k}{\mathrm{Im} \varphi_k(z)}\right)^{\gamma_k}+h_1\prod_{k=1}^n \left(\frac{\mathrm{Im} z_k}{\mathrm{Im} \psi_k(z)}\right)^{\gamma_k}\right)\right) \\
	\leqslant&C_{\gamma,\Omega}\left(\left(h_1 S_{\overline{\Omega},f^j}^{\gamma}+h_2 S_{\mathcal{H}^n \setminus\Omega,f^j}^{\gamma}\right)R_{\varphi,\psi}\right. \\
	&+\left.\rho(z)\left(h_2\prod_{k=1}^n \left(\frac{\mathrm{Im} z_k}{\mathrm{Im} \varphi_k(z)}\right)^{\gamma_k}+h_1\prod_{k=1}^n \left(\frac{\mathrm{Im} z_k}{\mathrm{Im} \psi_k(z)}\right)^{\gamma_k}\right)\right).
\end{align*}
Take $h_1=1,h_2=0$, we get
$$\left|f^j(\varphi(z))-f^j(\psi(z))\right|\prod_{k=1}^n (\mathrm{Im} z_k)^{\gamma_k}	\leqslant C_{\gamma,\Omega}\left(S_{\overline{\Omega},f^j}^{\gamma}R_{\varphi,\psi}+\rho(z)\prod_{k=1}^n \left(\frac{\mathrm{Im} z_k}{\mathrm{Im} \psi_k(z)}\right)^{\gamma_k}\right).$$

Use $J$ in the first case, then
$$\left|f^j(\varphi(z))-f^j(\psi(z))\right|\prod_{k=1}^n (\mathrm{Im} z_k)^{\gamma_k}	<C_{\gamma,\Omega}\left(\frac{\varepsilon}{2}+\frac{\varepsilon}{2}\right)=C_{\gamma,\Omega}\varepsilon.$$

\noindent $\bf{(c)}$ If $\varphi(z)\in\mathcal{H}^n \setminus\overline{\Omega},\psi(z)\in\Omega$, As in the previous case, we have
\begin{align*}
	&\left|f^j(\varphi(z))-f^j(\psi(z))\right|\prod_{k=1}^n (\mathrm{Im} z_k)^{\gamma_k}
	 \\ 
	& \leqslant C_{\gamma,\Omega}\left(\left(h_1 S_{\mathcal{H}^n \setminus\Omega,f^j}^{\gamma}+h_2 S_{\overline{\Omega},f^j}^{\gamma}\right)R_{\varphi,\psi}\right. \\
	&+\left.\rho(z)\left(h_2\prod_{k=1}^n \left(\frac{\mathrm{Im} z_k}{\mathrm{Im} \varphi_k(z)}\right)^{\gamma_k}+h_1\prod_{k=1}^n \left(\frac{\mathrm{Im} z_k}{\mathrm{Im} \psi_k(z)}\right)^{\gamma_k}\right)\right).
\end{align*}
Take $h_1=0,h_2=1$ to get
$$\left|f^j(\varphi(z))-f^j(\psi(z))\right|\prod_{k=1}^n (\mathrm{Im} z_k)^{\gamma_k}	\leqslant C_{\gamma,\Omega}\left(S_{\overline{\Omega},f^j}^{\gamma}R_{\varphi,\psi}+\rho(z)\prod_{k=1}^n \left(\frac{\mathrm{Im} z_k}{\mathrm{Im} \varphi_k(z)}\right)^{\gamma_k}\right)$$
and the subsequent steps proceed exactly as in the previous case.

\noindent $\bf{(d)}$ If $\varphi(z)\in\mathcal{H}^n \setminus\overline{\Omega},\psi(z)\in\mathcal{H}^n \setminus\overline{\Omega}$,
\begin{align*}
	&\left|f^j(\varphi(z))-f^j(\psi(z))\right|\prod_{k=1}^n (\mathrm{Im} z_k)^{\gamma_k} \\
	\leqslant& C_{\gamma,\Omega}\left(\left(h_1 S_{\mathcal{H}^n \setminus\Omega,f^j}^{\gamma}+h_2 S_{\mathcal{H}^n \setminus\Omega,f^j}^{\gamma}\right)\rho(z)\left(\prod_{k=1}^n \left(\frac{\mathrm{Im} z_k}{\mathrm{Im} \varphi_k(z)}\right)^{\gamma_k}+\prod_{k=1}^n \left(\frac{\mathrm{Im} z_k}{\mathrm{Im} \psi_k(z)}\right)^{\gamma_k}\right)\right. \\
	&+\left.S_{\mathcal{H}^n \setminus\Omega,f^j}^{\gamma}\rho(z)\left(h_2\prod_{k=1}^n \left(\frac{\mathrm{Im} z_k}{\mathrm{Im} \varphi_k(z)}\right)^{\gamma_k}+h_1\prod_{k=1}^n \left(\frac{\mathrm{Im} z_k}{\mathrm{Im} \psi_k(z)}\right)^{\gamma_k}\right)\right) \\
	\leqslant& 2C_{\gamma,\Omega}\| f^j \|_{A^{-\gamma}(\mathcal{H}^n)}\rho(z)\left(\prod_{k=1}^n \left(\frac{\mathrm{Im} z_k}{\mathrm{Im} \varphi_k(z)}\right)^{\gamma_k}+\prod_{k=1}^n \left(\frac{\mathrm{Im} z_k}{\mathrm{Im} \psi_k(z)}\right)^{\gamma_k}\right) \\
	\leqslant& 2C_{\gamma,\Omega}\rho(z)\left(\prod_{k=1}^n \left(\frac{\mathrm{Im} z_k}{\mathrm{Im} \varphi_k(z)}\right)^{\gamma_k}+\prod_{k=1}^n \left(\frac{\mathrm{Im} z_k}{\mathrm{Im} \psi_k(z)}\right)^{\gamma_k}\right) \\
	<&2C_{\gamma,\Omega}\varepsilon.
\end{align*}

In summary, we obtain
$$\left|f^j(\varphi(z))-f^j(\psi(z))\right|\prod_{k=1}^n (\mathrm{Im} z_k)^{\gamma_k}<2C_{\gamma,\Omega}\varepsilon,$$
and by the arbitrariness of $\varepsilon$, the proof is complete.

\end{proof}

\end{document}